\documentclass[12pt]{article}
\usepackage{}
\usepackage{amssymb}
\usepackage{amsmath}
\usepackage{amsthm}
\usepackage{amsfonts}
\usepackage{mathrsfs}
\usepackage{color, dsfont}
\allowdisplaybreaks
\usepackage{bm}
\usepackage{verbatim}
\usepackage{txfonts, graphicx, nicefrac}
\allowdisplaybreaks

  \def\bz{\beta}
    \def\dz{\delta}
    
\def\gz{\gamma}  
 
     \def\oz{\omega}
\def\pz{\pi}

        \def\sz{\sigma}
\def\tz{\tau}        
\def\vz{\varepsilon}

 \def\ddz{\Delta}
    
\def\ggz{\Gamma}  \def\ooz{\Omega}
\def\ppz{\Pi}

   \def\bq{{\mathscr{B}}}
   
   \def\fq{{\mathscr{F}}}

\def\qd{\quad}
\def\qqd{\qquad}

\def\lt{\left}
\def\rt{\right}

\def\PP{\mathbb{P}}
\def\EE{\mathbb{E}}

\newcommand{\mathsym}[1]{{}}

\def\leq{\leqslant}
\def\geq{\geqslant}

\newtheorem{thm}{Theorem}[section]
\newtheorem{prop}[thm]{Proposition}
\newtheorem{lem}[thm]{Lemma}
\newtheorem{rem}[thm]{Remark}

\newtheorem{defn}[thm]{Definition}
\newtheorem{exm}[thm]{Example}
\newtheorem{ass}[thm]{Assumption}

\numberwithin{equation}{section} \allowdisplaybreaks[4]

\def\prf{\medskip \noindent {\bf Proof}. }
\def\deprf{\quad $\square$ \medskip}

\def\be{\begin{equation}}
\def\de{\end{equation}}

\def\dear{\end{eqnarray*}}
\def\lb{\label}

\def\dps{\displaystyle}

\def\den{\end{enumerate}}

\def\d{\mathrm{d}}

\def\Ho{\mathcal{H}}

\def\PP{\mathbb{P}}
\def\EE{\mathbb{E}}

\def\NN{\mathbb{N}}

\def\one{\mathds{1}}

\begin{document}
	\date{}
	\pagestyle{plain}
	\title{Optimal stopping time on discounted semi-Markov processes \footnote{{\bf Funding:} This work was partly supported by the National Natural Science Foundation of China (No. 11931018, 61773411, 11701588) and the Guangdong Basic and Applied Basic Research Foundation (No. 2020B1515310021)}}
	\author{Fang Chen \footnote{School of Mathematics, Sun Yat-Sen University, Guangzhou 510275, China. Email: chenf76@mail2. sysu.edu.cn}, Xianping Guo \footnote{School of Mathematics, Sun Yat-Sen University, Guangzhou 510275, China. Email: mcsgxp@mail. sysu.edu.cn}, Zhong-Wei Liao \footnote{corresponding author. South China Research Center for Applied Mathematics and Interdisciplinary Studies, South China Normal University, Guangzhou 510631, China. Email: zhwliao@m.scnu.edu.cn}}
	\date{}
	\maketitle \underline{}
	
	{\bf Abstract:} This paper attempts to study the optimal stopping time for semi-Markov processes (SMPs) under the discount optimization criteria with unbounded cost rates. In our work, we introduce an explicit construction of the equivalent semi-Markov decision processes (SMDPs). The equivalence is embodied in the value functions of SMPs and SMDPs, that is, every stopping time of SMPs can induce a policy of SMDPs such that the value functions are equal, and vice versa. The existence of the optimal stopping time of SMPs is proved by this equivalence relation. Next, we give the optimality equation of the value function and develop an effective iterative algorithm for computing it. Moreover, we show that the optimal and $\vz$-optimal stopping time can be characterized by the hitting time of the special sets. Finally, to illustrate the validity of our results, an example of a maintenance system is presented in the end.
	

\setlength{\baselineskip}{0.25in}

\section{Introduction}

Optimal stopping theory is an important branch of the intersection of probability and control theory, which aims to find the optimal stopping time of stochastic systems according to a certain criterion. The optimal stopping theory of Markov processes have been widely used in finance, such as the pricing of American options, see the monographs \cite{BR11, CRS91, PS06} and the references therein. About the discrete-time Markov processes, Dochviri \cite{D95} constructed the corresponding relation between inhomogeneous and homogeneous Markov processes and solved the value function and the $\vz$-optimal stopping time; Nikolaev \cite{N99} proposed the multi-objective stopping time problem on Markov chains and proved the existence of $\vz$-optimal stopping times. About the continuous-time Markov processes, Zhitlukhin \& Shiryaev \cite{ZS14} gave the existence conditions of the optimal stopping time for unbounded reward functions; B\"auerle \& Popp \cite{BP18} studied the stopping time problems under the risk-sensitive criteria and obtained the optimality equation of the value function and the explicit expression of optimal stopping times; Ye \cite{Y17} proved the existence of optimal stopping time and gave the formula for the corresponding value function under the discount criteria.

It is well known that the sojourn time of discrete-time Markov processes is constant, while that of continuous-time Markov processes satisfies the exponential distribution. However, in practical application, the sojourn time not satisfy either of these two situations. The purpose of semi-Markov processes (SMPs) is to relax the condition of sojourn time,
which means that it satisfies the general distribution. This paper investigates the optimal stopping time of SMPs. As is well known, SMPs are a kind of general dynamic programming models with applications in many areas, see \cite{BR11, CO10, HG10, HG11, JN06, LO01, LYZ15, ZS14} for instance. For the optimal stopping time of SMPs, to the best of our knowledge, there is only a relevant research, see \cite{BG93}. Boshuizen \& Gouweleeuw \cite{BG93} studied the optimal stopping time problem of SMPs under discount criterion and solved the optimal stopping time by the dynamic programming method. It is worth noting that most of methods of optimal stopping problems are based on the martingale methods given by Snell \cite{S52}. In contrast, according to the particularities of SMPs, we use the techniques, which are introduced by B\"auerle \& Rieder \cite{BR11} for discrete-time Markov processes, to study the optimal stopping time of SMPs. Intuitively, at each jump epoch of the SMPs, the system has two options, continue or stop. Therefore, we can construct an action set with only two points $A = \{0, 1\}$, which are corresponding to the actions of the decision maker, i.e. continue ($a = 0$) or stop ($a = 1$). Hence, the optimal stopping problems of SMPs are equivalent to the semi-Markov decision processes (SMDPs), and finding the optimal stopping time of SMPs is equivalent to finding the optimal policy of SMDPs. Using this equivalence and the results of SMDPs given by Huang \& Guo \cite{HG11}, we can solve the optimal stopping time of SMPs.

The main contributions of the present paper are as follows. (1). In our model, we do not need to assume that the cost functions are continuous, monotonic or convex. That is much weaker than the hypothesis given in \cite{BG93, S52, Y17}. The only assumption of our model is Assumption \ref{s2-ass1}, which is the standard regular condition in SMPs, see \cite{HG10, HG11, P94, R70}. (2). We give a novel and explicit construction of SMDPs, which is equivalent to the optimal stopping problem of SMPs. The equivalence is embodied in the value functions of SMPs and SMDPs, that is, every stopping time of SMPs can induce a policy of SMDPs such that the value functions are equal, and vice versa, (see Proposition \ref{s3-prop1} and Theorem \ref{s3-thm1}). (3). Besides the existence of optimal stopping time, the calculation of the value function deserves more attention from the practical application. Hence, we introduce the optimality equation of the value function and derive an effective iterative algorithm for computing it, see Theorem \ref{s3-thm2}. Furthermore, through the optimality equation, we prove that the optimal stopping time can be characterized by the hitting time of a special set. (4). In order to numerical calculation, we also give the concept of the $\vz$-optimal stopping time and a sufficient condition to ensure that the $\vz$-optimal stopping time is optimal, see Theorem \ref{s3-thm4}. The significance of this result is that we can replace the value function with a computable approximate function, such that the $\vz$-optimal stopping time is equal to the optimal one. 

The rest of our paper is organized as follows. In Section 2, we describe the optimal stopping problems of SMPs and then give the regular condition (Assumption \ref{s2-ass1}). In Section 3, we present the explicit constructions of SMDPs, which are equivalent to the optimal stopping problems of SMPs. The equivalences between stopping time of SMPs and policy of SMDPs are introduced in Proposition \ref{s3-prop1} and Theorem \ref{s3-thm1}.  In Section 4, according to these equivalences, the optimality equation and the iterative algorithm of the value function are given in Theorem \ref{s3-thm2}. Moreover, we obtain the explicit expression of the optimal stopping time, see Theorem \ref{s3-thm3}. To numerical calculation, we give the concept of $\vz$-optimal stopping times and prove that under some conditions, the $\vz$-optimal stopping time is equal to the optimal one of SMPs. Finally, we illustrate the validity of our results by an example of a maintenance system in Section 5.

\section{The optimal stopping time model}

This paper studies the model of SMPs on a denumerable state space $S$ with transition mechanism $Q (t, j | i)$. Here and in what follows, $Q (t, j | i)$ is always assumed that
\begin{itemize}
	\item[(i)] given any $t \in [0, \infty)$, $Q (t, \cdot | \cdot)$ is a sub-stochastic kernel on $S$ given $S$;
	\item[(ii)] given any $i, j \in S$, $Q(\cdot, j| i)$ is a non-decreasing right continuous real-valued function on $[0, \infty)$, which satisfies $Q (0, j| i) =0$;
	\item[(iii)] $P (\cdot | \cdot) := \lim_{t \to \infty} Q (t, \cdot | \cdot)$ is a stochastic kernel on $S$ given $S$.
\end{itemize}
The evolution of the model is given as follows. At the beginning $t_0 = 0$, the system occupies the state $i_0 \in S$. Subsequently, the system remains in $i_0$ for time $t_1$ and then jumps to $i_1 \in S$ governed by the kernel $Q (t_1, i_1 | i_0)$. To describe the history (or trajectory, pathway) of SMPs, we introduce the measurable space $(\ooz, \bq (\ooz))$, which is based on the Kitaev construction (see \cite{K86, KR95}),
$$
\ooz = \lt\{ (i_0, t_1, i_1, \ldots, t_n, i_n \ldots) : i_0, i_n \in S, t_n \in [0, \infty), n \geq 1 \rt\},
$$
and $\bq(\ooz)$ is the corresponding Borel $\sz$-algebra. The element $\oz \in \ooz$ is known as a pathway of the system. The histories of SMPs up to the $n$-th jump epoch are
\be\lb{hn}
h_0 = i_0, \qd h_{n+1} = (i_0, t_1, i_1, \ldots, t_{n+1}, i_{n+1}), \qd n \geq 0.
\de
Denote by $\Ho_n$ the set of all histories $h_n$ up to the $n$-th jump epoch, which is endowed with the Borel $\sz$-algebra $\bq (\Ho_n)$. For each $n\geq 0$, $\oz = (i_0, t_1, i_1, \ldots, t_n, i_n, \ldots) \in \ooz$, define
$$
X_n (\oz) = i_n, \qd T_0(\oz) = 0, \qd T_{n+1} (\oz) = t_{n+1}
$$
where $X_n$ denotes the state at the $n$-th jump epoch and $T_{n+1}$ denotes the sojourn time between the $n$-th jump epoch and $(n+1)$-th jump epoch. Denote by
\be\lb{H_n}
H_n(\omega):=(X_0(\omega),T_0(\omega),X_1(\omega),\ldots,T_n(\omega),X_n(\omega)).
\de
In what follow, the argument $\oz$ is always omitted except some special informational statements. For each $n \geq 0$, denote by $\fq_n = \sz (H_n)$ the natural fluid of SMPs. For each $i \in S$, by the Tulcea theorem (see \cite[Proposition C.10]{HL96}), there exists a unique probability measure $\PP_i$ on $(\ooz, \bq(\ooz))$ such that, for each $t \in [0, \infty)$, $j \in S$ and $h_n=(i_0,t_1,i_1,\ldots,t_n,i_n) \in \Ho_n$ $(n \geq 0)$ it holds that
\begin{align}
& \PP_i ( T_0 =0, X_0 =i ) =1, \lb{pp-1}\\
& \PP_i (T_{n+1} \leq t, X_{n+1}=j | H_n = h_n) = Q( t, j | i_n ). \lb{pp-2}
\end{align}
Hence, $(\ooz, \bq(\ooz), \PP_i)$ becomes the probability space of SMPs, which is equipped with the natural fluid $\{ \fq_n, n \geq 0 \}$. Denote by $\EE_i$ the expectation with respect to $\PP_i$. Given any $n\geq 0$, we define the $n$-th jump epoch as
$$
S_n = \sum_{k=0}^n T_k.
$$
To ensure the regularity of SMPs, we give the following assumption.
\begin{ass}\lb{s2-ass1}
	There exist constants $\delta > 0$ and $\epsilon > 0$, such that
	\begin{equation}\lb{ass}
	\sum_{j \in S} Q (\delta, j | i) \leq 1 - \epsilon, \qd \forall i \in S.
	\end{equation}
\end{ass}
\noindent The Assumption \ref{s2-ass1} is a standard regular condition widely used in SMPs and SMDPs, see \cite{HG10, HG11, P94, R70}, for instance. According to \cite{HG11}, the Assumption \ref{s2-ass1} implies that
$$
\PP_i ( \lim_{n \to \infty} S_n = \infty )=1, \qd \forall i \in S.
$$
Corresponding to the process $\{ (T_n, X_n), n \geq 0\}$, we define an underlying continuous-time state process $\{ X(t), t \in [0, \infty) \}$ by
\be\lb{Xt}
X (t) = X_n, \qd S_n \leq t < S_{n+1}.
\de
Refer to Limnios and Oprisan \cite{LO01} for more details about the constructions of $\{ X(t), t \in [0, \infty) \}$ and the properties given in (\ref{pp-1}) and (\ref{pp-2}).

Next step, we introduce the optimal stopping time problems of SMPs. A mapping $\tz: \ooz \to \mathbb{N} \cup \{+\infty\}$ is called a $\fq_n$-stopping time, if for each $n \geq 0$, it holds that
\be\lb{stopping}
\{ \oz \in \ooz: \tz (\oz) = n \} \in \fq_n.
\de
Denote by $\ggz$ the set of all $\fq_n$-stopping times. In the absence of ambiguity, we write $\fq_n$-stopping time as stopping time. Let $c (i)$ and $g(i)$ be the nonnegative real-valued functions on $S$, which represent the cost rate function and the terminal cost function respectively. Fixed any discount factor $\beta>0$, the infinite horizon discounted cost is defined as:
\begin{displaymath}
R_{\tau} :=
\left\{
\begin{aligned}
&{ \int_{0}^{S_{\tz}} e^{-\beta t} c(X(t)) \d t + g ( X(S_\tz) ) e^{-\beta S_{\tz}}, } & &\tz< {+\infty}; \\
&{ \int_{0}^{+\infty} e^{-\beta t} c(X(t)) \d t, } & & \tz = {+\infty}. \\
\end{aligned}
\right.
\end{displaymath}
The infinite horizon expected discounted cost of a stopping time $\tau \in \ggz$ is given by
\begin{equation}\lb{V}
V^\tau (i) := \EE_{i} [R_\tau],\qqd i\in S.
\end{equation}
\begin{defn}
The function $V^* (i) := \inf_{\tz \in \ggz} V^{\tau} (i)$ is called the value function (or minimum expected discounted cost) of SMPs. A stopping time $\tz^* \in \ggz$ is called optimal if it achieves the infimum, i.e. $V^{\tz^*} (i) = V^* (i) = \inf_{\tz \in \ggz} V^{\tau} (i)$, for all $i \in S$.
\end{defn}

The main purpose of this paper is to find an optimal stopping time and give an algorithm for computing the value function $V^* (i)$.

\section{The equivalent semi-Markov decision process}

In this section, we will introduce the equivalent SMDPs corresponding to the original optimal stopping problem of SMPs. Intuitively, in the SMPs, stop or continue can be considered as a special action in the corresponding SMDPs. This intuition gives us an idea to construct the SMDPs.

The details about the constructions of SMDPs are given as follows. Here and in what follow, we always use ``$\ \hat{\cdot} \ $'' to distinguish the corresponding SMDPs from the original SMPs. The model of SMDPs is introduced by the four-tuple:
\be\lb{SMDPs}
\lt\{ \hat{S}, (A (i) \subset A), \hat{Q} (t, j | i, a), \hat{c}(i, a) \rt\}.
\de
The state space $\hat{S} := S \cup \{ \ddz \}$ is a denumerable space including the space $S$ of SMPs and a virtual state $\ddz$. The action space $A(i)$, which denotes the set of admissible actions at state $i\in \hat{S}$, is defined as
	$$
A (i):=
	\left\{
	\begin{array}{ll}
	\{0, 1 \}, & i\in S; \\
	\{ 1 \}, & i= \ddz,
	\end{array}
	\right.
	$$
where the action $0$ means continuation and $1$ means stop. The action space $A = \cup_{i \in \hat{S}} A (i)=\{0,1\}$ is finite. Denote by $K := \{ (i, a) : i \in \hat{S}, a \in A(i) \}$ the set of feasible state-action pairs. The semi-Markov kernel $\hat{Q} (t, j | i, a)$ of the SMDPs is given by
\be\lb{Q}
\hat{Q} (t, j | i, a) :=
\left\{
\begin{array}{ll}
	Q(t, j | i), & i \in S, j \in S, a = 0; \\
	\one_{ \lt[1 , +\infty \rt) } (t), & i \in \hat{S}, j = \ddz, a=1; \\
	0, & \text{otherwise},
\end{array}
\right.
\de
where $Q(t, j | i)$ is the kernel of SMPs and $\one_{E}$ is the indicator function on the set $E$. Finally, the cost rate function of SMDPs is defined as
\be\lb{r}
\hat{c}(i, a) :=
\left\{
\begin{array}{ll}
	c (i), & i \in S, a = 0; \\
	\bz g(i) / (1- e^{- \bz}), & i \in S, a=1; \\
	0, & i =  \Delta, a =1,
\end{array}
\right.
\de
where $c(i)$ is the cost rate function and $g(i)$ is the terminal cost function of SMPs.

The definitions of the history of SMDPs and the history-dependent policy are exactly the same as in \cite{HG10,HG11}, but for the ease of reading, we repeat it here. The trajectory space of SMDPs is defined as
$$
\hat{\ooz} = \lt\{(i_0, a_0, t_1, i_1, a_1, \ldots, t_n, i_n, a_n, \ldots) : \text{$ t_{m+1}\in [0,+\infty)$ and $(i_m,a_m)\in K$ for $m \geq 0$}\rt\},
$$
which is equipped with the Borel $\sz$-algebra $\bq (\hat{\ooz})$. Moreover, The histories of SMDPs up to the $n$-th jump epoch have the form
\be\lb{h'n}
\hat{h}_0 = i_0, \qd \hat{h}_{n+1} = \lt(i_0, a_0, t_1, i_1, \ldots,a_{n},  t_{n+1}, i_{n+1} \rt), \qd n \geq 0.
\de
Denote by $\hat{\Ho}_n$ the set of all histories $\hat{h}_n$, which is endowed with the Borel $\sz$-algebra $\bq (\hat{\Ho}_n)$. For each $\hat{\oz} =(i_0, a_0, t_1, i_1, a_1, \ldots, t_n, i_n, a_n, \ldots)\in \hat{\ooz}$, let
$$
\hat{X}_n(\hat{\omega})=i_n,\  \hat{A}_n(\hat{\omega})=a_n, \ \hat{T}_0(\hat{\omega})=0,\  \hat{T}_{n+1}(\hat{\omega})=t_{n+1}, \ \hat{S}_n(\hat{\omega}) := \sum_{k=0}^n \hat{T}_k(\hat{\omega}), \qd \forall n\geq0,
$$
 and $\hat{H}_n = (\hat{X}_0, \hat{A}_0, \hat{T}_1, \hat{X}_1, \hat{A}_1, \ldots, \hat{T}_n, \hat{X}_n)$. Similar to (\ref{Xt}), we define the continuous-time processes $\hat{X}(t)$, $\hat{A}(t)$ and the $n$-th jump epoch as following
\be
\hat{X} (t) = \hat{X}_n, \qd \hat{A} (t) = \hat{A}_n, \qd \hat{S}_n \leq t < \hat{S}_{n+1}.
\de
The definition of a deterministic history-dependent policy of SMDPs is given below, which specifies a decision rule to select actions.

\begin{defn}
A deterministic history-dependent policy is a sequence $\pz=\{ f_n, n \geq 0 \}$ of measurable functions $f_n:\hat{\Ho}_n \to A$ satisfying
$$
f_n(\hat{h}_n)\in A(i_n),\qd \forall \,\hat{h}_n=(i_0,a_0,t_1,i_1,\ldots,a_{n-1},t_{n},i_n)\in \hat{\Ho}_n, \, n\geq 0.
$$
In particular, the deterministic history-dependent policy is called deterministic stationary if $f_n$ are independent of $n$. Write $\pz=\{f,f,\ldots\}$ as $f$ for simplicity. Denote by $\ppz_{DH}$ and $\ppz_{DS}$ the sets of all deterministic history-dependent and deterministic stationary policies, respectively.
\end{defn}

For each $i \in \hat{S}$, the expected discounted cost of a policy $\pz \in \ppz_{DH}$ is defined as
$$
U^\pz (i) = \hat{\EE}_i^\pz \lt[ \int_0^\infty e^{-\beta t} \hat{c} (\hat{X}(t), \hat{A}(t)) \d t \rt],
$$
where $\hat{\EE}_i^\pz$ is the expectation depended on the state $i$ and the policy $\pz$, which is guaranteed by the Tulcea theorem (see \cite[Proposition C.10]{HL96} or \cite[Section 2]{HG10}). The value function of SMDPs is given as $U^* (i) = \inf_{\pz \in \ppz_{DH}} U^\pz (i)$.

Next step, we will focus on the relationship between the stopping times of SMPs and the policies of SMDPs. For each history of SMPs until the $n$-th jump epoch $h_n \in \Ho_n$ given in (\ref{hn}), we define a map $M_n$ as
\be\lb{Mn}
M_n (h_n) = \lt(i_0, 0, t_1, i_1, 0, \ldots, t_n, i_n \rt) \in \hat{\Ho}_n.
\de
The action $``0"$ (means continuation) added to the equation (\ref{Mn}) indicates that the system has been running incessantly until the $n$-th jump epoch. Obviously, $M_n (\Ho_n) := \{ M_n(h_n) : h_n \in \Ho_n \} \in \bq( \hat{\Ho}_n)$. Generally, define $M_n (C) := \{ M_n (h_n) : h_n \in C \} \in \bq( \hat{\Ho}_n)$ for each subset $C \in \bq (\Ho_n)$. For each stopping time $\tz \in \ggz$, we can introduce a policy $\pz_{\tz}$ in the following way.

\begin{defn}\lb{s3-p-t}
	Given any stopping time $\tz \in \ggz$ and $n \geq 0$, let
\be\lb{b_n}
B_n^\tau:=\big\{H_n(\omega): \omega \in \ooz, \tz (\oz) = n \big\},
\de
where $H_n$ is given in (\ref{H_n}). For each history $\hat{h}_n = (i_0, a_0, t_1, i_1,\ldots, a_{n-1},  t_n, i_n) \in \hat{\Ho}_n$, define
	$$
	f_n^{\tz} (\hat{h}_n) :=
	\left\{
	\begin{array}{ll}
	\one_{B_n^{\tz}} (i_0, t_1, i_1, \ldots, t_n, i_n), & \hat{h}_n \in M_n (\Ho_n); \\
	1, & \hat{h}_n \in \hat{\Ho}_n \setminus M_n (\Ho_n).
	\end{array}
	\right.
	$$
The policy $\pz_\tz := \{ f_n^\tz, n \geq 0 \}$ is called the policy induced by the stopping time $\tz$.
\end{defn}

\begin{lem}\lb{s3-lem1}
	For each stopping time $\tz \in \ggz$, the induced policy $\pi_\tau$ is a deterministic history-dependent policy of the corresponding SMDPs.
\end{lem}

\prf By Definition \ref{s3-p-t}, it holds that $f_n^\tau (\hat{h}_n) \in A(i_n)$. Then, we just need to consider the measurability. Noting that $B_n^\tau \in \bq (\Ho_n)$, we have $\big\{ \hat{h}_n \in \hat{\Ho}_n : f_n^\tau (\hat{h}_n) = 0 \big\} = M_n (\Ho_n) \bigcap M_n \big( (B_n^\tau)^c \big) \in \bq (\hat{\Ho}_n)$ and $\big\{ \hat{h}_n \in \hat{\Ho}_n : f_n^\tau (\hat{h}_n) =1 \big\}=\hat{\Ho}_n\setminus \big\{ \hat{h}_n \in \hat{\Ho}_n : f_n^\tau (\hat{h}_n) = 0 \big\}\in \bq (\hat{\Ho}_n)$. Hence, the policy $\pz_\tz := \{ f_n^\tz, n \geq 0 \}$ becomes a deterministic history-dependent policy of the corresponding SMDPs.
\deprf

The Definition \ref{s3-p-t} and Lemma \ref{s3-lem1} say that given any stopping time of SMPs, we can construct a history-dependent policy of SMDPs.  Next, on the contrary, given any history-dependent policy of SMDPs, we construct a corresponding stopping time of SMPs, see Definition \ref{s3-t-p} and Lemma \ref{s3-lem2} below.

\begin{defn}\lb{s3-t-p}
	Given any deterministic history-dependent policy $\pi=\{f_n,n\geq 0\}$. For each $\oz \in \ooz$, we define
	$$
	\tz_\pz (\oz) := \inf \lt\{ n \in \mathbb{N} : f_n ( M_n (H_n (\oz) ) = 1 \rt\},
	$$
	where $\inf \{ \emptyset \} := +\infty$ and $M_n$ is given in (\ref{Mn}). Then $\tau_\pi$ is called the stopping time induced by the policy $\pi$.
\end{defn}

\begin{lem}\lb{s3-lem2}
	For each deterministic history-dependent policy $\pi=\{f_n,n\geq 0\}$ of SMDPs, the induced stopping time $\tau_\pi$ is a stopping time.	
\end{lem}

\prf Note that for each $n \geq 0$, the random variable $H_n = (X_0, T_1, X_1, \ldots,T_n X_n)$, the mapping $M_n$ and the function $f_n$ are measurable in their corresponding spaces. Hence, we have
$$
\{ \tz_\pz = n \} = \lt( \bigcap_{k=0}^{n-1} \{ f_k ( M_k (H_k) ) = 0 \} \rt) \bigcap \lt\{ f_n ( M_n (H_n) ) = 1 \rt\} \in \fq_n,
$$
which implies that $\tau_\pi$ is a $\fq_n$-stopping time.
\deprf

Those results give us an idea that the stopping times of SMPs and the policies of the SMDPs are one-to-one correspondence. Hence, we give the following proposition to verify this idea.

\begin{prop}\lb{s3-prop1}
	Given any stopping time $\tau \in \Gamma$ of SMPs, it holds that
	\be\lb{s3-prop1-1}
	\tau = \tau_{\pi_\tau},
	\de
	where $\pi_\tau = \{ f_n^\tau, n\geq 0\}$ is the policy induced by $\tz$ and $\tau_{\pi_\tau}$ is the stopping time induced by $\pi_\tau$.
\end{prop}

\prf Fix any $n \geq 0$. For each $\oz \in \ooz$, we have
\be\lb{s3-prop-2}
f_n^\tau (M_n (H_n (\oz))) =\one_{B_n^\tau} (H_n (\oz)) = \one_{\{\tau=n\}} (\omega).
\de
Hence, by Definition \ref{s3-t-p}, it holds that $\tau_{\pi_\tau} (\oz) = \inf \big\{ n\in \mathbb{N} : f_n^\tz (M_n (H_n (\oz))) = 1 \big\} = \inf \big\{ n\in \mathbb{N} : \one_{\{\tau=n\}} = 1 \big\}$. Thus, we obtain that
$$
\big\{\tau_{\pi_\tau}  =n \big\} = \big\{\omega\in \Omega : \one_{\{\tau = k \}} (\omega) =0, 0\leq k\leq n-1,\one_{\{\tau=n\}} (\omega) =1 \big\} = \big\{ \tau=n \big\}.
$$
By the arbitrariness of $n \geq 0$, we obtain $\tau=\tau_{\pi_\tau}$.
\deprf

In short, Proposition \ref{s3-prop1} ensure that the relationship between the stopping times of SMPs and the policies of the SMDPs is one to one correspondence. Moreover, for any stopping time $\tau$, the infinite horizon expected discounted cost of $\tau$ is equivalent to the expected discounted cost of the policy $\pi_\tau$ that corresponds to it, see Theorem \ref{s3-thm1} below.

\begin{thm}\lb{s3-thm1}
	Given any $\tau\in \Gamma$, let $\pi_\tau = \{f_n^\tau, n\geq 0\}$ be the induced policy of $\tz$. Then, we have
	\be\lb{s3-thm1-2}
	V^\tau(i)=U^{\pi_\tau}(i), \qd \forall i \in S.
	\de
\end{thm}

\prf For each $\hat{\oz} = (i_0, a_0, t_1,\ldots,  i_n, a_n,t_{n+1} \ldots) \in \hat{\ooz}$, $\hat{H}_n (\hat{\omega})= (i_0, a_0,t_1,i_1 \ldots,a_{n-1}, t_n, i_n)$ is the history of SMDPs up to the $n$-th jump epoch. Denote by $\{ C_n, n \geq 0\}$ and $C$ the subsets of $\hat{\ooz}$, which are defined as
\begin{align*}
C_n & := \lt\{ \hat{\oz} \in \hat{\ooz} : \inf \lt\{ k \in \NN: f_k^\tz (\hat{H}_k(\hat{\oz})) = 1 \rt\} =n \rt\}, \qd n \geq 0;\\
C & := \lt\{ \hat{\oz} \in \hat{\ooz} : \text{ $f_k^\tz (\hat{H}_k(\hat{\oz})) = 0$ for all $k \geq 0$}\rt\}.
\end{align*}
It is easy to know that $\{C, C_n, n\geq 0 \}$ is a partition of $\hat{\ooz}$. Then, using the monotone convergence theorem, we obtain
\be\lb{pf-thm1-1}
U^{\pi_\tau}(i)=\sum_{k=0}^{+\infty} \hat{\EE}_i^{\pi_\tau} \lt[\one_{C_k} \int_0^\infty e^{- \beta t} \hat{c}(\hat{X}(t), \hat{A}(t)) \d t \rt] +\hat{\EE}_i^{\pi_\tau} \lt[ \one_{C} \int_0^\infty e^{- \beta t} \hat{c}(\hat{X}(t), \hat{A}(t)) \d t \rt].
\de

To calculate the first item of (\ref{pf-thm1-1}), for each $k\geq 0$, we use
\begin{align}\lb{pf-thm1-2}
&\hat{\EE}_i^{\pi_\tau} \lt[\one_{C_k} \int_0^\infty e^{- \beta t} \hat{c}(\hat{X}(t), \hat{A}(t)) \d t \rt] = \hat{\EE}_i^{\pi_\tau} \lt[\one_{C_k} \sum_{m=0}^{+ \infty} \int_{\hat{S}_m}^{\hat{S}_{m+1}} e^{- \beta t} \hat{c}(\hat{X}(t), \hat{A}(t)) \d t \rt] \notag \\
& = \sum_{m=0}^{+ \infty} \frac{1}{\beta} \hat{\EE}_i^{\pi_\tau} \lt[\one_{C_k} \lt(e^{-\beta \hat{S}_m} - e^{-\beta\hat{S}_{m+1}} \rt) \hat{c}(\hat{X}_m, \hat{A}_m) \rt] = \sum_{m=0}^{k} \frac{1}{\beta} \hat{\EE}_i^{\pi_\tau} \lt[\one_{C_k} \lt(e^{-\beta \hat{S}_m} - e^{-\beta \hat{S}_{m+1}} \rt) \hat{c}(\hat{X}_m, \hat{A}_m) \rt].
\end{align}
The third equality of (\ref{pf-thm1-2}) is based on the definition of $\hat{Q}$ given in (\ref{Q}). In fact, the semi-Markov kernel $\hat{Q}$ satisfies $\hat{Q}(t, \Delta | i,1) =1$ for each $t \geq 1$ and $i \in \hat{S}$, and the action set satisfies $A(\ddz) = \{ 1\}$. These mean that once the action $1$ is selected, the process will jump to state $\ddz$ with probability one after a unit time, and then stay in $\ddz$ forever, i.e.
\be
\hat{A}_k(\hat{\omega})=1, \  \hat{A}_l(\hat{\omega})=1, \  \hat{X}_l(\hat{\omega})=\Delta, \qd \text{$\hat{\omega} \in C_k$ and $l \geq k+1$}.
\de
Since $\hat{c}(\Delta,1)=0$, it holds that $\one_{C_k} \hat{c}(\hat{X}_l,\hat{A}_l) = 0$ for each $l \geq k+1$. In next step, our goal is to show that
\be\lb{pf-thm1-3}
\hat{\EE}_i^{\pi_\tau} \lt[ \one_{C_k} \lt( e^{-\beta \hat{S}_{m}} - e^{-\beta \hat{S}_{m+1}} \rt) \hat{c} (\hat{X}_m, \hat{A}_m) \rt] =
\left\{
\begin{array}{ll}
	\EE_i \lt[ \one_{\{\tau=m\}} \lt( e^{-\beta S_{m}} - e^{-\beta S_{m+1}} \rt) c (X_m) \rt], & 0 \leq m < k; \\
	\EE_i \lt[ \beta \one_{\{\tau=k\}} \lt( e^{-\beta S_k} \rt) g (X_k) \rt], & m = k.
\end{array}
\right.
\de
In the beginning, let's discuss the relationship between $C_k$ and $\{ \tz = k \}$. By the definition of $C_k$,  $\hat{\oz} \in C_k $ if and only if $f_m^\tau ( \hat{H}_m (\hat{\oz})) = 0$ ($0 \leq m < k$) and $f_k^\tau ( \hat{H}_k (\hat{\oz})) =1$. Conversely, according to Definition \ref{s3-t-p} and Proposition \ref{s3-prop1} , we have
	$$
	\tz (\oz) = \tz_{\pz_{\tz}} (\oz) = \inf \lt\{n \in \NN: f_n^\tz (M_n (H_n (\oz))) = 1 \rt\},
	$$
which means that $\one_{\{\tau=k\}}=\prod_{m=0}^{k-1}(1-f_m^\tz (M_m (H_m )))\times f_k^\tz (M_k (H_k ))$. For each $0 \leq m < k$, we have
\begin{align}\lb{pf-thm1-4}
& \hat{\EE}_i^{\pi_\tau} \lt[ \one_{C_k}  \lt( e^{-\alpha \hat{S}_{m}} - e^{-\alpha \hat{S}_{m+1}} \rt) \hat{c}(\hat{X}_m, \hat{A}_m) \rt] \notag \\
& =\sum_{i_0\in S}\delta_{i}(i_0)\sum_{i_1\in \hat{S}}\int_0^\infty \hat{Q}(\d t_1,i_1|i_0,0)\cdots  \sum_{i_k \in \hat{S}}\int_0^\infty \hat{Q}(\d t_k,i_k|i_{k-1},0)\hat{c}(i_m,0) \notag \\
&\qd \times \lt[ e^{ - \beta \sum_{n=1}^{m} t_n} - e^{-\beta \sum_{n=1}^{m+1} t_n} \rt] \prod_{n=0}^{k-1}\lt(1-f_n^\tau(i_0,0,t_1, \ldots, 0, t_n, i_n) \rt) \times f_k^\tau (i_0, 0, t_1, \ldots, 0, t_k, i_k),
\end{align}
where $\sum_{n=1}^0 t_n:=0$. Then, by the definitions of $\hat{Q}$, $r$, and $f_n^\tau$ given in (\ref{Q}), (\ref{r}), and Definition \ref{s3-p-t}, we obtain
\begin{align}\lb{pf-thm1-5}
& \hat{\EE}_i^{\pi_\tau} \lt[ \one_{C_k}  \lt( e^{-\beta \hat{S}_{m}} - e^{-\beta \hat{S}_{m+1}} \rt) \hat{c}(\hat{X}_m, \hat{A}_m) \rt] \notag \\
& = \sum_{i_0\in S}\delta_{i}(i_0)\sum_{i_1 \in S} \int_0^{\infty} Q (\d t_1, i_1 | i_0 ) \cdots \sum_{i_k \in S} \int_0^{\infty} Q( \d t_k, i_k | i_{k-1}) \lt[ e^{(- \beta \sum_{n=1}^{m} t_n)} - e^{(-\beta \sum_{n=1}^{m+1} t_n)} \rt] \notag \\
&\qd \times c (i_m) \prod_{n=0}^{k-1}\lt(1-f_n^\tau(M_n(i_0,t_1,i_1,\ldots,t_n,i_n))\rt) \times f_k^\tau (M_k(i_0,t_1,i_1,\ldots,t_k,i_k))\notag\\
& = \EE_i \lt[ \one_{\{\tau=k\}} c (X_m) \lt(e^{-\beta S_m} - e^{-\beta S_{m+1}} \rt) \rt].
\end{align}
In the same way, we can calculate (\ref{pf-thm1-3}) in the case $m=k$. The only thing to be careful about is that $\hat{Q}(t, \Delta | i, 1) = \one_{ \lt[ 1, +\infty \rt) } (t)$ for all $i \in \hat{S}$. It means that the system will occupy the state $i$ with a unit time if the action $1$ is selected at the state $i$. Since $\hat{c} (i, 1) = \bz g(i) / (1 - e^{-\bz})$, we have
\begin{align}\lb{pf-thm1-6}
& \hat{\EE}_i^{\pi_\tau} \lt[ \one_{C_k} \lt( e^{-\bz \hat{S}_{k}} - e^{-\bz \hat{S}_{k+1}} \rt) \hat{c} (\hat{X}_k, \hat{A}_k) \rt] \notag \\
& = \sum_{i_0\in \hat{S}}\delta_{i}(i_0) \sum_{i_1 \in \hat{S}} \int_0^{\infty} \hat{Q} (\d t_1, i_1 | i_0 , 0) \cdots \sum_{i_k \in \hat{S}} \int_0^{\infty} \hat{Q} (\d t_k, i_k | i_{k-1}, 0) \frac{\bz g(i_k)}{1 - e^{-\bz}} \notag \\
& \qd \times e^{( -\bz \sum_{n=1}^k t_n)} \lt( 1 - e^{-\bz} \rt) \prod_{n=0}^{k-1} \lt (1-f_n^\tau (i_0, 0, t_1, \ldots, 0, t_n, i_n ) \rt) f_k^\tau (i_0, 0, t_1, \ldots, 0, t_k, i_k) \notag \\
& = \bz \sum_{i_0 \in S} \delta_{i} (i_0) \sum_{i_1 \in S} \int_0^{\infty} Q(\d t_1, i_1 | i_0) \cdots \sum_{i_k \in S} \int_0^{\infty} Q (\d t_k, i_k | i_{k-1}) e^{( -\bz \sum_{n=0}^k t_n )} g(i_k) \notag \\
&\qd \times \prod_{n=1}^{k-1}\lt(1-f_n^\tau(M_n(i_0,t_1, \ldots, t_n, i_n) ) \rt) \times f_k^\tau (M_k(i_0, t_1, \ldots, t_k, i_k) ) \notag \\
& = \EE_i \lt[ \bz \one_{\{\tau=k\}} e^{- \bz S_k} g(X_k) \rt].
\end{align}
Hence, (\ref{pf-thm1-5}) and (\ref{pf-thm1-6}) imply that (\ref{pf-thm1-3}) holds. Moreover, by (\ref{pf-thm1-2}) and (\ref{pf-thm1-3}), we have
\begin{align}\lb{pf-thm1-7}
& \hat{\EE}_i^{\pi_\tau} \lt[\one_{C_k} \int_0^\infty e^{- \beta  t} \hat{c}(\hat{X}(t), \hat{A}(t)) \d t \rt] = \sum_{m=0}^k \frac{1}{\beta } \hat{\EE}_i^{\pz_{\tz}} \lt[ \one_{C_k} \lt(e^{-\beta \hat{S}_m} - e^{-\beta  \hat{S}_{m+1}} \rt) \hat{c} (\hat{X}_m, \hat{A}_m) \rt] \notag \\
& = \sum_{m=0}^{k-1} \frac{1}{\beta } \EE_i \lt[\one_{\{ \tz = k \}} \lt(e^{-\beta  S_m} - e^{-\beta S_{m+1}} \rt) c (X_m) \rt] + \frac{1}{\beta } \EE_i \lt[ \beta \one_{\{ \tz = k \}} e^{-\beta S_k} g(X_k) \rt] \notag \\
& = \EE_i \lt[ \one_{\{ \tz = k\}} R_{\tz} \rt].
\end{align}
Next, we calculate the second item of (\ref{pf-thm1-1}). Noting that $\{C, C_n, n\geq 0 \}$ is a partition of $\hat{\ooz}$, for all $k\geq 0$ we obtain that
$$
\one_C = \lt( 1-\sum_{n=0}^{+\infty } \one_{C_n} \rt) \times \prod_{m=0}^{k}(1-\one_{C_m}) = \prod_{m=0}^{k}(1 - \one_{C_m}) - \sum_{n=k+1}^{+\infty } \one_{C_n},
$$
which implies that
\begin{align}\lb{pf-thm1-8}
 \hat{\EE}_i^{\pz_{\tz}} \bigg[ \one_C \int_0^\infty\hat{c} (\hat{X}(t), \hat{A}(t)) \d t \bigg]
= &\sum_{k=0}^{+ \infty} \frac{1}{\beta} \bigg\{
 \hat{\EE}_i^{\pi_\tau}  \bigg[\prod_{m=0}^{k}(1-\one_{C_m}) \hat{c}(\hat{X}_k,\hat{A}_k)(e^{-\beta \hat{S}_k}-e^{-\beta \hat{S}_{k+1}})  \bigg] \notag \\
&- \sum_{n=k+1}^{+\infty} \hat{\EE}_i^{\pi_\tau}  \bigg[\one_{C_n} \hat{c}(\hat{X}_k,\hat{A}_k)(e^{-\beta\hat{S}_k}-e^{-\beta\hat{S}_{k+1}})  \bigg]  \bigg\}
\end{align}
Similar to (\ref{pf-thm1-8}), we decompose $\EE_i \bigg[ \one_{\{\tau=+\infty\}}R_\tau  \bigg]$ as
\begin{align}\lb{pf-thm1-9}
\EE_i \bigg[ \one_{\{\tau=+\infty\}}R_\tau  \bigg] & = \sum_{k=0}^{+ \infty} \frac{1}{\beta} \Bigg\{ \EE_i \bigg[\prod_{m=0}^{k}(1-\one_{\{\tau=m\}})c(X_k)(e^{-\beta S_k}-e^{-\beta S_{k+1}}) \bigg] \notag \\
& \qd - \sum_{n=k+1}^{+\infty}\EE_i \bigg[\one_{\{\tau=n\}} c(X_k)(e^{-\beta S_k}-e^{-\beta S_{k+1}})  \bigg]  \Bigg\}.
\end{align}
For the first item of (\ref{pf-thm1-9}), similar to (\ref{pf-thm1-5}), it holds that
$$
\EE_i \bigg[\prod_{m=0}^{k}(1-\one_{\{\tau=m\}})c(X_k)(e^{-\beta S_k}-e^{-\beta S_{k+1}})\bigg]
= \hat{\EE}_i^{\pi_\tau}  \bigg[\prod_{m=0}^{k}(1-\one_{C_m}) \hat{c}(\hat{X}_k,\hat{A}_k)(e^{-\beta \hat{S}_k}-e^{-\beta\hat{S}_{k+1}})  \bigg].
$$
Hence, together with (\ref{pf-thm1-5}), (\ref{pf-thm1-8}), and (\ref{pf-thm1-9}), it holds that
\be\lb{pf-thm1-10}
\EE_i \bigg[ \one_{\{\tau=+\infty\}}R_\tau  \bigg] =\hat{\EE}_i^{\pz_{\tz}} \bigg[ \one_C \int_0^\infty e^{-\beta t}\hat{c} (\hat{X}(t), \hat{A}(t)) \d t \bigg] .
\de
Finally, according to (\ref{pf-thm1-1}), (\ref{pf-thm1-7}) and (\ref{pf-thm1-10}), we obtain
\begin{align*}
& U^{\pi_\tau}(i)  = \sum_{k=0}^{+\infty} \EE_i \lt[ \one_{\{ \tz = k\}} R_{\tz} \rt] + \EE_i \lt[ \one_{\{\tz = \infty\}} R_{\tz} \rt] = V^\tz (i).
\end{align*}
The proof of this theorem is completed.
\deprf

\section{The iterative algorithm and optimal stopping times}

Those results given in Section 3 state that for each stopping time $\tau$ of SMPs, there exists an equivalent policy $\pi_\tau$ of the corresponding SMDPs. Hence, we can analyze the value function $V^*$ and the optimal stopping time $\tz^*$ of SMPs through the conclusions of the corresponding SMDPs. To do so, we give some results about the value function $U^*$ and the optimal policy of SMDPs by using the results in \cite{HG10}. Note that the regular condition (Assumption \ref{s2-ass1}) is needed.

\begin{lem}\lb{u-thm}
Suppose that Assumption \ref{s2-ass1} holds. For the corresponding SMDPs defined in (\ref{SMDPs}), the following statements hold.
	
(a). Let $U^*_{-1}(i) \equiv 0$. For each $n \geq 0$ and  $i \in \hat{S}$, define a sequence of functions $U_n^*$ on $\hat{S}$ as
\begin{equation}\lb{U-function}
U_{n}^*(i)=\min_{a\in A(i)}\lt\{  \hat{c}(i,a)\int_0^\infty e^{-\beta t} \lt(1-\sum_{j\in \hat{S}}\hat{Q}(t,j|i,a) \rt) \d t + \sum_{j\in \hat{S}} \int_0^\infty e^{-\beta t}\hat{Q}(\d t,j|i,a)U_{n-1}^*(j) \rt\}.
\end{equation}
Then, for each $i \in \hat{S}$, $U_{n}^*(i)$ are non-descending and the value function of SMDPs $U^* (i)$ satisfies $U^*(i)=\lim_{n \rightarrow \infty}U^*_n(i)$ for each $i \in \hat{S}$. Moreover, it holds that
\begin{equation}\lb{U-function-2}
U^*(i)=\min_{a\in A(i)}\lt\{  \hat{c}(i,a)\int_0^\infty e^{-\beta t} \lt( 1-\sum_{j\in \hat{S}}\hat{Q}(t,j|i,a) \rt)\d t + \sum_{j\in \hat{S}}\int_0^\infty e^{-\beta t}\hat{Q}(\d t,j|i,a)U^*(j)\rt\}.
\end{equation}

(b). A deterministic stationary policy $f \in \Pi_{DS}$ is optimal if and only if for each $i\in \hat{S}$
\begin{equation}\lb{f-op}
U^*(i)=  \hat{c}(i,f(i))\int_0^\infty e^{-\beta t} \lt( 1-\sum_{j\in \hat{S}}\hat{Q}(t,j|i,f(i)) \rt)\d t +
\sum_{j\in \hat{S}}\int_0^\infty e^{-\beta t}\hat{Q}(\d t,j|i,f(i))U^*(j).
\end{equation}

(c). There exists a deterministic stationary policy $f^*\in \Pi_{DS}$ satisfying (\ref{f-op}), which means $f^*$ is optimal.
\end{lem}

\prf According to Assumption \ref{s2-ass1}, there exists $\delta>0$ and $\epsilon>0$ such that for all $i\in S$, $\sum_{j\in S}Q(\delta, j|i)\leq 1-\epsilon$. Let $\delta^*:=\min\{ \delta, 1/2 \}$. By the definition of $\hat{Q}$ in (\ref{Q}), we have that
\begin{equation}\lb{pf-u}
\sum_{j\in \hat{S}}\hat{Q} (\delta^*, j | i, a) =
\left\{
\begin{array}{ll}
\sum_{j\in S} Q(\delta^*, j | i) \leq 	\sum_{j\in S} Q(\delta, j | i) \leq  1-\epsilon , & i \in S,  a = 0; \\
\one_{ \lt[ 1, +\infty \rt) } (\delta^*) =0  \leq 1-\epsilon, & i \in \hat{S},  a=1,
\end{array}
\right.
\end{equation}
which implies $\sum_{j\in \hat{S}}\hat{Q}(\delta^*,j|i,a)\leq  1-\epsilon$ for all $(i,a) \in K$. Hence, the regular condition of \cite{HG10} holds. Part $(a)$, $(b)$ and $(c)$ come from \cite[Theorem 4.1]{HG10}, \cite[Theorem 3.1]{HG10} and \cite[Theorem 3.2]{HG10} respectively.
\deprf

With these preparation in hand, we can give our main results on SMPs. That is, the existence of the optimal stopping time of SMPs and then we give an iterative algorithm for computing the value function $V^*$, see Theorem \ref{s3-thm2} below. To establish the algorithm, we define an operator $T$ from $\mathbb{M}$ to $\mathbb{M}$, where $\mathbb{M}$ denotes the set of all non-negative real-valued functions on $S$, i.e. for each  $V\in \mathbb{M}$, $T V(i)$ is defined by
	\be\lb{it-al}
	T V (i) := \min \lt\{ g(i), c(i) \int_0^{\infty} e^{-\bz t} \lt(1 - \sum_{j\in S} Q(t, j | i) \rt) \d t + \sum_{ j \in S}\int_0^\infty e^{-\bz t} Q(\d t, j | i ) V (j) \rt\}.
	\de
From this expression, we can verify that $T$ is a monotone operator, i.e. $T V_1 \leq T V_2$ if $V_1 \leq V_2$. 
\begin{thm}\lb{s3-thm2}
Suppose that Assumption \ref{s2-ass1} holds. Then, the following statements hold.
\begin{itemize}
	\item[(a).] There exists an optimal stopping time of SMPs.
	\item[(b).] Let $V^*_{-1} (i) \equiv 0$ and $V^*_{n} (i) = T V^*_{n-1} (i)$ for each $n \geq 0$ and $i \in S$. Then, for each $i \in S$, $V_{n}^*(i)$ are non-descending and the value function of SMPs $V^*$ satisfies $V^*(i)=\lim_{n \rightarrow \infty}V^*_n(i)$ and is the solution to the following optimality equation
	$$
	V^*(i) = T V^*(i), \qd \forall i \in S.
	$$
\end{itemize}
\end{thm}

\prf $(a)$. Under Assumption \ref{s2-ass1}, Lemma \ref{u-thm} says that there exists a policy $f^*\in \Pi_{DS} \subset \Pi_{DH}$ such that for all $i\in S$, $U^* (i) = U^{f^*} (i)$. Using Proposition \ref{s3-prop1} and Theorem \ref{s3-thm1}, we have
\begin{equation}\lb{pf-s3-thm2-1}
U^*(i) = \inf_{\pz \in \Pi_{DH}} U^{\pz} (i) = \inf_{\tau \in \Gamma} V^{\tz} (i) \leq V^{\tau_{f^*}}(i) = U^{f^*} (i) = U^*(i), \qd \forall i\in S,
\end{equation}
where $\tau_{f^*}$ is the stopping time induced by the policy $f^*$. Hence, we have
$$
V^*(i) = \inf_{\tau \in \Gamma} V^{\tz} (i) = V^{\tau_{f^*}}(i) , \qd \forall i\in S,
$$
which means that $\tau_{f^*}$ is the optimal stopping time of SMPs.

$(b)$. When $i = \ddz \in \hat{S}$, the functions $U^*_n (i)$ given in Lemma \ref{u-thm} satisfy that
$$
U_{n}^*(\Delta) = r(\Delta,1) \int_0^\infty e^{-\beta t} \lt( 1-\hat{Q}(t,\Delta|\Delta,1) \rt) \d t +
\int_0^\infty e^{-\beta t} \hat{Q}(\d t,\Delta | \Delta, 1) U^*_{n-1} (\Delta), \qd n \geq 0.
$$
Once the initial value $U^*_{-1}(\Delta)=0$ is given, we have $U^*_n(\Delta)=0$ for each $n \geq 0$. Hence, we only consider the case $i \in S$. By Lemma \ref{u-thm} and (\ref{pf-s3-thm2-1}), we have
\be\lb{u-eq-v}
V^*(i)=U^*(i)=\lim_{n \rightarrow +\infty}U_{n}^*(i), \qd \forall i\in S.
\de
Next, by induction, we aim to prove that
\be\lb{pf-s3-thm2}
U_n^*(i) = V^*_n(i), \qd \text{$\forall i\in S$ and $n \geq -1$}.
\de
Clearly, (\ref{pf-s3-thm2}) holds for $n = -1$. Assume that it holds for $n = k-1$, we now consider the case $n=k$. Note that $A(i) = \{0, 1\}$ for each $i \in S$, the expression (\ref{U-function}) implies that
\begin{align}\lb{pf-s3-thm2-3}
U^*_{k} (i) = &\min_{a \in A(i)} \lt\{\hat{c}(i,a) \int_0^\infty e^{-\beta t} (1- \sum_{j \in \hat{S}} \hat{Q}(t,j | i,a) ) \d t + \sum_{j\in \hat{S}} \int_0^\infty e^{-\beta t} \hat{Q} (\d t,j | i,a) U_{k-1}^*(j) \rt\} \notag \\
=& \min \Bigg\{ \hat{c}(i,0) \int_0^\infty e^{-\beta t} \lt(1 - \sum_{j \in \hat{S}} \hat{Q}(t,j | i,0) \rt) \d t + \sum_{j \in \hat{S}} \int_0^\infty e^{-\beta t} \hat{Q} (\d t,j | i,0) U_{k-1}^*(j) , \notag \\
& \qd \qd \hat{c}(i,1) \int_0^\infty e^{-\beta t} \lt( 1 - \sum_{j\in \hat{S}} \hat{Q}(t,j | i,1) \rt) \d t +\sum_{j \in \hat{S}} \int_0^\infty e^{-\beta t} \hat{Q} (\d t,j | i,1) U_{k-1}^*(j) \Bigg\}.
\end{align}
For the second item of (\ref{pf-s3-thm2-3}), by the definitions of $\hat{Q}$ in (\ref{Q}) and $r$ in (\ref{r}), we have
\begin{align*}
& \hat{c}(i,1) \int_0^\infty e^{-\beta t} \lt( 1 - \sum_{j\in \hat{S}} \hat{Q}(t,j | i,1) \rt) \d t + \sum_{j \in \hat{S}} \int_0^\infty e^{-\beta t} \hat{Q} (\d t,j | i,1) U_{k-1}^*(j) \\
& \qd = \hat{c}(i,1) \int_0^{\infty} e^{-\beta t} \lt( 1- \hat{Q}(t, \Delta| i,1) \rt) \d t + \int_0^\infty e^{-\beta t} \hat{Q} (\d t, \Delta | i,1) U^*_{k-1}(\Delta) \\
& \qd = \frac{\bz g(i)}{1 - e^{-\bz}} \int_0^{1} e^{-\beta t} \d t = g(i).
\end{align*}
Hence,
\begin{align*}
U^*_{k} (i) & = \min \lt\{\hat{c}(i,0) \int_0^\infty e^{-\beta t} \lt(1 - \sum_{j \in \hat{S}} \hat{Q}(t,j | i,0) \rt) \d t + \sum_{j \in \hat{S}} \int_0^\infty e^{-\beta t} \hat{Q} (\d t,j | i,0) U_{k-1}^*(j) , g(i) \rt\} \\
& = \min \lt\{ c(i) \int_0^{\infty} e^{-\beta t} \lt( 1- \sum_{j\in S} Q( t, j| i) \rt) \d t + \sum_{j \in S} \int_0^\infty e^{-\bz t} Q (\d t, j | i) V^*_{k-1} (j), g(i) \rt\} \\
& = T V^*_{k-1} (i) = V^*_{k} (i).
\end{align*}
According to (\ref{u-eq-v}) and (\ref{pf-s3-thm2}), the part $(b)$ of this theorem holds, and we complete the proof.
\deprf

From the perspective of practical application, it is not enough to only have the existence of the optimal stopping time of SMPs. We are more looking forward to giving a computable characterization of the optimal stopping time. In the following theorem, we give an optimal stopping time, which is equivalent to the hitting time of a special set.

\begin{thm}\lb{s3-thm3}
	Suppose that Assumption \ref{s2-ass1} holds. Let $V^* (i)$ be the value function of SMPs. Define a subset of $S$ by
	\be\lb{set-s}
	S^* := \lt\{ i \in S : g(i) =TV^*(i) \rt\}.
	\de
	Then for each $\omega = (i_0, t_1\ldots, i_n,t_{n+1} \ldots) \in \Omega$, the optimal stopping time is equal to the hitting time of $S^*$, that is
	\be\lb{tz*}
	\tau^*(\omega) =
	\left\{
	\begin{aligned}
	&+ \infty, & & \text{$S^* = \varnothing$}; \\
	&\inf \{ n \in \NN: i_n \in S^* \}, & & \text{$S^* \neq \varnothing$}. \\
	\end{aligned}
	\right.
	\de
\end{thm}

\prf According to Theorem \ref{s3-thm2}, we have
$$
V^*(i) =
\left\{
\begin{aligned}
&g(i),& &i \in  S^*; \\
&c(i) \int_0^{\infty} e^{-\beta t} \lt( 1- \sum_{j\in S}Q( t, j| i) \rt) \d t+ \sum_{ j \in S} \int_0^\infty e^{-\beta t} Q (\d t, j | i) V^*(j), & & i \in S \setminus S^*. \\
\end{aligned}
\right.
$$
Let $U^* (i)$ be the value function of the corresponding SMDPs, and its existence is guaranteed by Lemma \ref{u-thm}. Moreover, by (\ref{u-eq-v}), it holds that $V^* (i) = U^* (i)$ for each $i \in S$ and $U^*(\ddz) = 0$.

For each $i \in \hat{S}$, denote by $f^*(i) = \one_{S^* \cup \{\Delta\}} (i)$ a deterministic stationary policy of SMDPs. What we need to do is to verify that $f^* \in \Pi_{DS}$ is an optimal policy of SMDPs, that is, $f^*$ satisfies (\ref{f-op}) of Lemma \ref{u-thm}. Since $A(\Delta)=\{1\}$, $f^*(\Delta)=1$ satisfies (\ref{f-op}). For each $i \in S^*$, we have $f^*(i) =1$. Hence, the definitions of $\hat{Q}$ and $r$ imply that
\begin{align*}
&\hat{c}(i,f^*(i))\int_0^\infty e^{-\beta t} \lt( 1-\sum_{j\in \hat{S}}\hat{Q}(t,j|i,f^*(i)) \rt)\d t + \sum_{j\in \hat{S}} \int_0^\infty e^{-\beta t}\hat{Q} (\d t,j | i,f^*(i)) U^*(j) \\
& =g(i) = V^*(i) =U^*(i).
\end{align*}
For each $i \in S \setminus S^*$, we have $f^*(i) =0$. Again, using the definitions of $\hat{Q}$ and $r$, we have
\begin{align*}
& \hat{c}(i, f^*(i)) \int_0^\infty e^{-\beta t} \lt(1-\sum_{j\in \hat{S}}\hat{Q}(t,j|i,f^*(i)) \rt) \d t + \sum_{j\in \hat{S}} \int_0^\infty e^{-\beta t} \hat{Q}(\d t,j |i, f^*(i)) U^*(j) \\
&=c(i) \int_0^\infty e^{-\beta t} \lt( 1-\sum_{j\in S}Q(t,j|i) \rt) \d t + \sum_{ j \in S} \int_0^\infty e^{-\beta t} Q (\d t, j | i) V^*(j) = U^*(i).
\end{align*}
Hence, (\ref{f-op}) holds and $f^*$ is an optimal deterministic stationary policy of SMDPs. Denote by $\tz^* := \tau_{f^*}$ the stopping time induced by the policy $f^*$, which is an optimal stopping time of SMPs according to Theorem \ref{s3-thm2}. If $S^* \neq \varnothing$, according to Theorem \ref{s3-thm2}, for each $\omega=(i_0, t_1, i_1, \ldots, t_n, i_n,\ldots) \in \Omega$, it holds that
$$
\tau^* (\omega) = \tau_{f^*} (\omega) = \inf \{n \in \mathbb{N} : f^* ( i_n ) = 1\} = \inf \{n \in \mathbb{N} : i_n\in S^* \}.
$$
If $S^* = \varnothing$, for all $\oz \in \ooz$, $f^* (i_n) = \one_{\ddz} (i_n) = 0$ ($\forall n \geq 0$). Then,
$$
\tau^*(\omega) = \tau_{f^*} (\omega) = \inf \{n \in \mathbb{N} : f^* ( i_n ) = 1\} = + \infty,
$$
and the proof is achieved.
\deprf

The condition of Theorem \ref{s3-thm3} requires the value function $V^* (i)$ of SMPs, but in practical application, the value function is often unknown. Intuitively, we can replace the value function $V^* (i)$ by the approximation function $V_n^* (i)$, which is obtained by the iterative algorithm given in Theorem \ref{s3-thm2}. Therefore, the concept of optimal stopping time will be replaced by $\varepsilon$-optimal, that is the following definition.
\begin{defn}
Given any $\varepsilon>0$, a stopping time of SMPs $\tz^\varepsilon$ is called $\varepsilon$-optimal if it holds that $V^{\tz^\varepsilon} (i) - V^*(i) \leq \varepsilon$, for all $i \in S$, where $V^*(i)$ is the value function of SMPs and $V^{\tz^\varepsilon} (i)$ is the expected discounted cost given in (\ref{V}).
\end{defn}
The following theorem shows that for any $\vz >0$, we can iterate enough times and get a $\vz$-optimal stopping time under some conditions. For the convenience of statement, we give two notations, i.e. $||f|| := \sup_{x\in E} |f(x)|$ for any function $f$ defined on the set $E$; $\lfloor x \rfloor := \max \{ n \in \mathbb{N} : n \leq x \}$ for any $x\in [0,+\infty)$. 

\begin{thm}\lb{s3-thm4}
Suppose that Assumption \ref{s2-ass1} holds and $c(i)$ and $g(i)$ are bounded. For any $\varepsilon>0$, the number of iterations $N_\varepsilon$ is given by
\be\lb{N-epsilon}
N_{\varepsilon}:=\lt\lfloor \frac{\log(\varepsilon(\epsilon-e^{-\beta \delta^*}))-\log(\beta^{-1}||c||+||g||+1)}{\log(1-\epsilon+\epsilon e^{-\beta \delta^*})} \rt\rfloor,
\de
where $\epsilon$, $\delta$ are introduced in Assumption \ref{s2-ass1} and $\dz^* = \min\{\dz, 1/2\}$. Let $V^*_{N_\varepsilon} (i)$ be the $N_\varepsilon$-th step iterative function given in Theorem \ref{s3-thm2} and $S^\varepsilon$ be the subset of $S$ given by $S^\varepsilon := \lt\{ i \in S : g(i) = T V^*_{N_\varepsilon} (i) \rt\}$. Then, the following statements hold.
\begin{itemize}
\item[(a).] Denote by $\tau^\varepsilon$ the hitting time of $S_\varepsilon$, i.e.
\begin{equation*}
\tau^\varepsilon(\omega) =
\left\{
\begin{aligned}
	&+ \infty, & & \text{$S^\varepsilon = \varnothing$}; \\
	&\inf \{ n \in \NN: i_n \in S^\varepsilon \}, & & \text{$ S^\varepsilon  \neq \varnothing$}. \\
\end{aligned}
\right.
\end{equation*}
Then, $\tau^\varepsilon$ is an $\varepsilon$-optimal stopping time of SMPs.
\item[(b).] If it holds that
\be\lb{st-op}
\inf_{i \in S \setminus S^\varepsilon} \lt( g(i)-TV^*_{N_\varepsilon}(i) \rt) > \varepsilon,
\de
then, $\tau^\varepsilon$ is also the optimal stopping time of SMPs.
\end{itemize}
\end{thm}

\prf $(a)$. Using Lemma \ref{u-thm}, for each $i \in \hat{S}$ and $a \in A(i)$, we have
\begin{align}\lb{s3-thm4-1}
& \sum_{j\in \hat{S}} \int_0^\infty e^{-\beta t} \hat{Q} (t,j | i,a) \d t = \int_0^{\dz^*} e^{- \bz t} \lt( \sum_{j\in \hat{S}} \hat{Q} (\d t,j | i,a) \rt) + \int_{\dz^*}^\infty e^{- \bz t} \lt( \sum_{j\in \hat{S}} \hat{Q} (\d t,j | i,a) \rt) \notag \\
& \qd \leq \lt( 1- e^{- \bz \dz^*} \rt) \lt( \sum_{j\in \hat{S}} \hat{Q} (\dz^*, j | i,a) \rt) + e^{- \bz \dz^*} \leq 1-\epsilon+ \epsilon e^{-\beta \delta^*} < 1,
\end{align}
where the second inequality depends on (\ref{pf-u}). Denote by $\gamma := 1-\epsilon+ \epsilon e^{-\beta \delta^*} < 1$. In the same way, we obtain $\sum_{ j \in S} \int_0^\infty e^{-\bz t} Q(\d t, j | i ) \leq \gz$. Since $V_{n}^* (i) - V_{n-1}^* (i) \leq \lt\| V_{n}^* - V_{n-1}^* \rt\|$ for each $n \geq 0$ and $i \in S$,  by the definition of $T$ in (\ref{it-al}), we have
\begin{align*}
T V_{n}^* (i) & \leq \min \Bigg\{g(i), c(i) \int_0^{\infty} e^{-\bz t} \lt(1 - \sum_{j\in S} Q(t, j | i) \rt) \d t + \sum_{ j \in S} \int_0^\infty e^{-\bz t} V^*_{n-1} (j) Q(\d t, j | i ) \\
& \qqd \qqd \qqd + \lt\| V_{n}^* - V_{n-1}^* \rt\| \sum_{ j \in S} \int_0^\infty e^{-\bz t} Q(\d t, j | i ) \Bigg\} \\
& \leq \min \Bigg\{g(i), c(i) \int_0^{\infty} e^{-\bz t} \lt(1 - \sum_{j\in S} Q(t, j | i) \rt) \d t + \sum_{ j \in S} \int_0^\infty e^{-\bz t} V^*_{n-1} (j) Q(\d t, j | i ) \Bigg\} \\
& \qqd \qqd \qqd + \lt\| V_{n}^* - V_{n-1}^* \rt\| \sum_{ j \in S} \int_0^\infty e^{-\bz t} Q(\d t, j | i ) \\
& \leq T V_{n-1}^* (i) + \gz \lt\| V_{n}^* - V_{n-1}^* \rt\|.
\end{align*}
Let $U_n^* (i)$ be the function given in (\ref{U-function}) on Lemma \ref{u-thm}, and then using (\ref{pf-s3-thm2}), we have $U^*_n (i) = V_n^* (i)$ for all $i \in S$ and $n \geq -1$. Hence,
$$
U_{n}^* (i) - U^*_{n-1} (i) = T V_{n-1}^* (i) - T V_{n-2}^* (i) \leq \gz \lt\| V_{n-1}^* - V_{n-2}^* \rt\| = \gz \lt\| U_{n-1}^* - U_{n-2}^* \rt\|,
$$
which implies that (since $U^*_{-1} (i) \equiv 0$)
\be\lb{s3-thm4-2}
 \lt\| U^*_{n}- U^*_{n-1} \rt\| \leq \gamma^{n} \lt\| U^*_{0} \rt\| \leq \gamma^{n} \lt(\beta^{-1} \| c \| + \| g \|+1 \rt).
\de

For each $i \in \hat{S}$, we define $f^\varepsilon(i):=\one_{S^\varepsilon \cup \{\Delta\}}$. Similar to the proof of Theorem \ref{s3-thm3}, we obtain $f^\varepsilon \in \Pi_{DS}$. Moreover, for each $n \geq 1$, by induction, we show that
\begin{align}\lb{s3-thm4-5}
U^*_{N_\varepsilon+1} (i) & \geq \sum_{m=0}^{n-1}  \hat{\EE}^{f^\varepsilon}_i \lt[ (e^{-\beta \hat{S}_m} - e^{-\beta \hat{S}_{m+1}}) \hat{c} (\hat{X}_m, \hat{A}_m) \rt] + \hat{\EE}^{f^\varepsilon}_i \lt[ e^{-\beta \hat{S}_{n}} U^*_{N_\varepsilon+1}(\hat{X}_{n}) \rt] \notag \\
& \qd - \sum_{m=0}^{n-1} \gamma^{m+1} \lt\| U^*_{N_\varepsilon+1} - U^*_{N_\varepsilon} \rt\|.
\end{align}
The key of (\ref{s3-thm4-5}) is that, by Lemma \ref{u-thm}, it holds
\begin{align*}
& U^*_{N_\varepsilon + 1} (i_n) \geq \hat{c} (i_n, f^\varepsilon (i_n)) \int_0^\infty e^{-\beta t} \lt( 1-\sum_{i_{n+1} \in \hat{S}} \hat{Q} (t, i_{n+1} | i_n , f^\varepsilon(i_n)) \rt) \d t \\
& \qqd \qqd + \sum_{i_{n+1} \in \hat{S}} \int_0^\infty e^{-\beta t} \hat{Q} (\d t, i_{n+1} | i_n, f^\varepsilon(i_n)) U^*_{N_\varepsilon + 1} (i_{n+1}) - \gamma \lt\| U^*_{N_\varepsilon+1}-U^*_{N_\varepsilon} \rt\| \\
&= \hat{\EE}^{f^\varepsilon}_i \lt[ \lt(1-e^{-\beta (\hat{S}_{n+1} - \hat{S}_n)} \rt) \hat{c}(\hat{X}_n,\hat{A}_n) \Big| \hat{X}_n = i_n \rt] \\
& \qqd \qqd + \hat{\EE}^{f^\varepsilon}_i \lt[ e^{-\beta (\hat{S}_{n+1} - \hat{S}_n) } U^*_{N_\varepsilon+1}(\hat{X}_{n+1}) \Big| \hat{X}_n = i_n \rt] - \gamma \lt\| U^*_{N_\varepsilon+1}-U^*_{N_\varepsilon} \rt\|.
\end{align*}
Hence, the second item of (\ref{s3-thm4-5}) satisfies
\begin{align*}
\hat{\EE}^{f^\varepsilon}_i \lt[ e^{-\beta \hat{S}_{n}} U^*_{N_\varepsilon+1}(\hat{X}_{n}) \rt] & \geq \hat{\EE}^{f^\varepsilon}_i \lt[ \lt(e^{-\beta \hat{S}_{n}} -e^{-\beta \hat{S}_{n+1}} \rt) \hat{c} (\hat{X}_n, \hat{A}_n) \rt] \\
& \qd+ \hat{\EE}^{f^\varepsilon}_i \lt[ e^{-\beta \hat{S}_{n+1}} U^*_{N_\varepsilon+1}(\hat{X}_{n+1}) \rt] + \gz^{n+1} \lt\| U^*_{N_\varepsilon+1}-U^*_{N_\varepsilon} \rt\|.
\end{align*}
Hence, passing the limit $n \rightarrow \infty$, it holds that
\be\lb{s3-thm4-6}
U^*(i) \geq U^*_{N_\varepsilon+1} (i) \geq U^{f^\varepsilon}(i)- \frac{\gamma}{1-\gamma}||U^*_{N_\varepsilon+1}-U^*_{N_\varepsilon}||,  \qd \forall \ i\in \hat{S}.
\de
Using (\ref{s3-thm4-2}) and the definition of $N_{\vz}$, it can be verified that $\dps\frac{\gz}{1-\gz} \lt\| U^*_{N_\varepsilon+1} - U^*_{N_\varepsilon} \rt\| \leq \vz$. Then, for each $i \in \hat{S}$, we have $U^*(i) \geq U^{f^\varepsilon}(i) - \varepsilon$. In the same method of Theorem \ref{s3-thm3}, we can verify that the hitting time $\tau^\varepsilon$ of $S^{\vz}$ satisfies $\tau^\varepsilon = \tau_{f^\varepsilon}$, where $\tau_{f^\varepsilon}$ is the stopping time induced by $f^{\vz}$. Moreover, using Theorem \ref{s3-thm1} and Theorem \ref{s3-thm2}, we have
$$
V^*(i) = U^*(i) \geq U^{f^\varepsilon} (i) - \varepsilon = \EE_i \lt[ R_{\tau^\varepsilon} \rt] -\varepsilon = V^{\tau^\varepsilon} (i) - \vz,
$$
which means that $\tau^\varepsilon$ is an $\varepsilon$-optimal stopping time of SMPs.

$(b)$. If $S \setminus S^\varepsilon = \varnothing$, then $S = S^\varepsilon$ and the condition (\ref{st-op}) holds naturally. Hence, for each $i \in S$, we have $g(i) = T V^*_{N_\varepsilon} (i) \leq T V^* (i)$. That means $S = S^* = S^\varepsilon$, where $S^*$ given in (\ref{set-s}). 

Next, we consider the case $S \setminus S^\varepsilon \neq \varnothing$. Again, using the monotonicity of $T$, we have $S^\vz \subset S^*$. Conversely, by the definition of $V^*_{N_\varepsilon} (i)$ and (\ref{s3-thm4-6}), it holds that
$$
T V^*_{N_\varepsilon} (i) = V^*_{N_\varepsilon + 1} (i) \geq U^{f^\vz} (i) - \vz \geq U^* (i) - \vz = V^* (i) - \vz, 
$$
i.e. $V^* (i) \leq T V^*_{N_\varepsilon} (i) + \vz$ for each $i \in S$. For each $i \in S \setminus S^\varepsilon$, the condition (\ref{st-op}) implies that
$$
g(i) - T V^*(i) = g(i) - V^*(i) \geq g(i) - T V^*_{N_\varepsilon} - \varepsilon >0,
$$
which means $i \in S \setminus S^*$. Hence, we have $S^* \subset S^\varepsilon$ and then $S^* = S^\varepsilon$. Finally, we have $\tau^\varepsilon = \tau^*$, which is an optimal stopping time given in Theorem \ref{s3-thm3}. 
\deprf

\section{An application to a maintenance system}

In this section, we consider a specific example of the optimal stopping time of SMPs, that is, the maintenance system. We will illustrate how to calculate the value function and the optimal stopping time by the algorithm.

\begin{exm}
The repairable maintenance system is made up by three states, say $1$, $2$ and $3$, which represent ``normal operation'', ``minor failure'' and ``serious failure'' respectively. At each state $i \in \{1, 2, 3 \}$, the decision maker has two choices, either to maintain the system or to stop using it. If the decision maker chooses to maintain the system, it will incur maintenance cost with rate $c(i)$. After a random period of time, which obeys the exponential distribution with parameter $u (i) > 0$, the system transfers to state $j \in \{1, 2, 3 \}$ with probability $p_{ij}$. Otherwise, if the decision maker chooses to stop using the system, a terminal cost $g(i)$ will have to be paid. And then the system stops running and there is no need to pay any more.
	
We now model this problem as an optimal stopping time problem of SMPs. Let $S = \{1, 2, 3\}$ be the state space. The semi-Markov kernel is given by $Q(t , j | i) = p_{ij} \lt( 1 - e^{-u(i) t} \rt)$ for each $t \geq 0$ and $i, j \in S$, where $p_{ij}$ is the transition probability and $u(i)$ is the parameter of exponential distribution. The cost $c(i)$ and the terminal cost $g(i)$ are determined by the data in the system. Next, to illustrate the effectiveness of Theorem \ref{s3-thm2} and Theorem \ref{s3-thm3}, we consider a numerical example. Fixed the discount factor with $\bz = 0.05$. Assume that $p_{ij}$, $u(i)$, $c(i)$ and $g(i)$ satisfy that
	\begin{align*}
	p_{ij} & =
	\left(
	\begin{array}{ccc}
	0.8 & 0.15 & 0.05 \\
	0.6 & 0.2 & 0.2 \\
	0.1 & 0.1 & 0.8
	\end{array}
	\right); \qd u(i) = \lt(0.1, 2, 1 \rt); \\
	c(i) &= (5, 30, 80); \qd g(i) = (300, 350, 400).
	\end{align*}
Under the data above, the Assumption \ref{s2-ass1} holds in this example with $\delta=1$ and $\epsilon=\min\{e^{-0.1},e^{-1},e^{-2}\}$. Hence, we have $\gamma=1-e^{-2}+e^{-2.025}$. Then, let $V_{-1}^* (i) =0$ and $V^*_{n} (i)= T V^*_{n-1} (i)$ for $n \geq 0$, where $T$ is the operator defined in (\ref{it-al}). Denote by $V^* (i)=\lim_{n \to \infty} V_n^* (i)$. When the number of iterations $n=84$ in {\rm Matlab}, we get that $\|V_{85}^*-V_{84}^*\|\leq 10^{-12}$. Let $\varepsilon=10^{-8}$. Then we have 
	\begin{equation*}
	V^*_{84} (1) = 147.6923, \qd V^*_{84} (2) = 222.5641,  \qd V^*_{84}(3) = 400;
	\end{equation*}
and $\|V^* - V^*_{84} \| \leq \vz$. Next step, we consider the set $S^\varepsilon = \lt\{ i \in S : g(i) =TV^*_{84}(i) \rt\}$ given in Theorem \ref{s3-thm4}. To do so, by numerical calculation, we have
\begin{align*}
g(1)-c(1) \int_0^{\infty} e^{-0.05 t} \lt( 1- \sum_{j\in S}Q( t, j| 1) \rt) \d t+ \sum_{ j \in S} \int_0^\infty e^{-0.05 t} Q (\d t, j | 1) V^*_{84}(j)&>10^{-8} , \\
g(2)-c(2) \int_0^{\infty} e^{-0.05 t} \lt( 1- \sum_{j\in S}Q( t, j| 2) \rt) \d t + \sum_{ j \in S} \int_0^\infty e^{-0.05 t} Q (\d t, j | 1) V^*_{84}(j)&>10^{-8} , \\
c(3) \int_0^{\infty} e^{-0.05 t} \lt( 1- \sum_{j\in S}Q( t, j| 3) \rt) \d t + \sum_{ j \in S} \int_0^\infty e^{-0.05 t} Q (\d t, j | 1) V^*_{84}(j) = 416 &> g(3).
\end{align*}
Hence, $S^*=S^\varepsilon = \{ 3\}$ and the optimal stopping time of this maintenance system is given by
$$
\tau^*(\omega) = \inf\{n \in \NN : i_n = 3 \}, \qd \omega = (i_0, t_1\ldots, i_n, t_{n+1},\ldots) \in \Omega.
$$
\end{exm}

\bibliographystyle{alpha}

\end{document}